\documentclass{article}    
\usepackage{amstext}
\def\qed{$\rlap{$\sqcap$}\sqcup$}
\usepackage{latexsym}

\begin{document}           

{\ }\\
\begin{center}
{\huge {\bf Level algebras of type 2}} \\ [.250in]
{\large FABRIZIO ZANELLO\\
Dipartimento di Matematica, Universit\`a di Genova, Genova, Italy.\\E-mail: zanello@dima.unige.it}
\end{center}
{\large

{\ }\\
\\
ABSTRACT. In this paper we study standard graded artinian level algebras, in particular those whose socle-vector has type 2. Our main results are: the characterization of the level $h$-vectors of the form $(1,r,...,r,2)$ for $r\leq 4$; the characterization of the minimal free resolutions associated to each of the $h$-vectors above when $r=3$; a sharp upper-bound (under certain mild hypotheses) for the level $h$-vectors $(1,r,...,a,2)$ of arbitrary codimension $r$ and type 2, which depends on the next to last entry $a$.\\
\\
\section{Introduction}
\indent

In this paper we study level algebras, focusing our attention mainly on those of type 2. The algebras $A=R/I$ are standard graded artinian quotients of the polynomial ring $R=k[x_1,...,x_r]$, where $k$ is a field of characteristic zero.\\\indent
Before introducing our results, let us recall briefly the main definitions we need. The {\it $h$-vector} of $A$ is $h(A)=h=(h_0,h_1,...,h_e)$, where $h_i=\dim_k A_i$ and $e$ is the last index such that $\dim_k A_e>0$. Since we may suppose that $I$ does not contain non-zero forms of degree 1, $r=h_1$ is defined as the {\it codimension} of $A$.\\\indent 
The {\it socle} of $A$ is the annihilator of the maximal homogeneous ideal $\overline{m}=(\overline{x_1},...,\overline{x_r})\subseteq A$, namely $soc(A)=\lbrace a\in A {\ } \mid {\ } a\overline{m}=0\rbrace $. Since $soc(A)$ is a homogeneous ideal, we define the {\it socle-vector} of $A$ as $s(A)=s=(s_0,s_1,...,s_e)$, where $s_i=\dim_k soc(A)_i$. Note that $h_0=1$, $s_0=0$ and $s_e=h_e>0$. The integer $e$ is called the {\it socle degree} of $A$ (or of $h$). The {\it type} of the socle-vector $s$ (or of the algebra $A$) is type$(s)=\sum_{i=0}^es_i$.\\\indent 
If $s=(0,0,...,0,s_e=t)$, we say that the algebra $A$ is {\it level}. In particular, if $t=1$, $A$ is {\it Gorenstein}. With a slight abuse of notation, we will sometimes refer to an $h$-vector as Gorenstein (or level) if it is the $h$-vector of a Gorenstein (or level) algebra.\\
\\\indent
Gorenstein $h$-vectors of codimension $r=2$ and $r=3$ have been characterized by Stanley in $[St]$ (see also $[Ma]$ for $r=2$), while in higher codimensions only a few results on the possible Gorenstein $h$-vectors are known (see $[IS]$ for $r=4$ and $[BI]$, $[Bo]$ and $[BL]$ for larger $r$'s. See also $[Ha]$, $[MN]$ and $[CI2]$). Also, in $[Di]$, Diesel characterized the minimal free resolutions (MFR's, in brief) associated to the Gorenstein $h$-vectors of codimension 3.\\\indent
In general, level $h$-vectors were first studied by Stanley (see $[St2]$). Iarrobino, in his 1984 paper $[Ia2]$, determined all the level $h$-vectors of codimension 2. See $[Ia2]$ and $[FL]$ for the $h$-vectors of (level) {\it compressed} algebras, i.e. algebras having the (entry by entry) maximal possible $h$-vector, given codimension and socle-vector. (See also this author's works $[Za]$ and $[Za2]$.)\\\indent
Recently, level $h$-vectors have received a lot of attention: see $[CI]$, $[BG]$, $[GHS]$, $[Ia3]$, $[GHMS]$ and $[Ia]$.\\
\\\indent
Our paper is structured as follows: in Section 2 we prove a stronger form of a decomposition theorem shown by Geramita {\it et al.} in $[GHMS]$ and, above all, we characterize a class of level $h$-vectors of type 2 and codimension less than or equal to 4. Our result is the first characterization of a class of non-Gorenstein $h$-vectors of codimension $r>2$, and concerns one of the \lq \lq simplest" cases: that of level $h$-vectors $h$ having type 2 and next to last entry equal to $r$, where $r\leq 4$. These $h$-vectors have a very nice description: $h$ is level if and only if $h$ can be written as the sum of a Gorenstein $h$-vector of codimension $r-1$ and $(0,1,1,...,1)$.\\\indent 
In Section 3 we investigate the structure of the level $h$-vectors of type 2, say $(1,r,...,a,2)$, by a different technique, especially studying the intervals where the entries of such $h$-vectors may range. Our main result is a sharp upper-bound for these $h$-vectors, which (of course) depends on the next to last entry $a$ (under the only restrictions $a\geq r$, and $r\leq 5$ if $a=r$).\\\indent 
Finally, in Section 4 we determine the possible minimal free resolutions associated to the class of level $h$-vectors of the form $h=(1,3,...,3,2)$ (which we have characterized in Section 2). In particular, it turns out that there is a unique MFR associated to $h$, unless $h=(1,3,3,...,3,3,2)$, when there are exactly two MFR's. Our characterization of the possible MFR's for the $h$-vectors above is the first result of this kind for a class of non-Gorenstein $h$-vectors of codimension $r>2$.\\\indent 
The results obtained in this paper are part of the author's Ph.D. dissertation, written at Queen's University (Kingston, Ontario, Canada), under the supervision of Professor A.V. Geramita.\\
\\\indent 
Recall now the main facts of the theory of Inverse Systems which we will use throughout the paper. For a complete introduction, we refer the reader to $[Ge]$ and $[IK]$.\\\indent 
Let $S=k[y_1,...,y_r]$, and consider $S$ as a graded $R$-module where the action of $x_i$ on $S$ is partial differentiation with respect to $y_i$.\\\indent 
There is a one-to-one correspondence between artinian algebras $R/I$ and finitely generated $R$-submodules $M$ of $S$, where $I=M^{-1}$ is the annihilator of $M$ in $R$ and, conversely, $M$ is the $R$-submodule of $S$ which is annihilated by $I$ (cf. $[Ge]$, Remark 1), p. 17).\\\indent 
If $R/I$ has socle-vector $s$, then $M$ is minimally generated by $s_i$ elements of degree $i$, for $i=1,...,e$, and the $h$-vector of $R/I$ is given by the number of linearly independent derivatives in each degree obtained by differentiating the generators of $M$ (cf. $[Ge]$, Remark 2), p. 17).\\\indent 
In particular, level algebras of type $t$ and socle degree $e$ correspond to $R$-submodules of $S$ minimally generated by $t$ elements of degree $e$.\\ 
\\\indent 
The next theorem is a well-known result of Macaulay, which we will often use in the sequel.\\
\\\indent 
{\bf Definition-Remark 1.1.} Let $n$ and $i$ be positive integers. The {\it i-binomial expansion of n} is $$n_{(i)}={n_i\choose i}+{n_{i-1}\choose i-1}+...+{n_j\choose j},$$ where $n_i>n_{i-1}>...>n_j\geq j\geq 1$.\\\indent 
Under these hypotheses, the $i$-binomial expansion of $n$ is unique (e.g., see $[BH]$, Lemma 4.2.6).\\\indent 
Furthermore, define $$n^{<i>}={n_i+1\choose i+1}+{n_{i-1}+1\choose i-1+1}+...+{n_j+1\choose j+1}.$$
\\\indent 
{\bf Theorem 1.2} (Macaulay). {\it Let $h=(h_i)_{i\geq 0}$ be a sequence of non-negative integers, such that $h_0=1$, $h_1=r$ and $h_i=0$ for $i>e$. Then $h$ is the $h$-vector of some standard graded artinian algebra if and only if, for every $d$, $1\leq d\leq e-1$, $$h_{d+1}\leq h_d^{<d>}.$$}
\\\indent 
{\bf Proof.} See $[BH]$, Theorem 4.2.10. (This theorem holds, with appropriate modifications, for any standard graded algebra, not necessarily artinian.){\ }{\ }\qed \\
\\
\indent A sequence of non-negative integers which satisfies the growth condition of Macaulay's theorem is called an {\it $O$-sequence}.\\
\\
\section{A class of level $h$-vectors of type 2}
\indent

The following is a decomposition theorem due to Geramita {\it et al.}. We present it in a form which is the most useful to our purposes.\\
\\
\indent {\bf Theorem 2.1} ($[GHMS]$). {\it Let $F_1,...,F_t\in S$ be any linearly independent forms of degree $e$, and let $M=\langle F_1,...,F_t \rangle $ and $N=\langle F_1,...,F_{t-1}\rangle $ be two submodules of $S$. Then the reverse of the difference between the $h$-vectors of $A=R/M^{-1}$ and $R/N^{-1}$ is an $O$-sequence, which is the $h$-vector of a quotient of $A$.\\}
\\
{\bf Proof.} See $[GHMS]$, Lemma 2.8 and Theorem 2.10.{\ }{\ }\qed \\
\\\indent 
As an immediate consequence we have:\\
\\
\indent {\bf Corollary 2.2} ($[GHMS]$). {\it Let $h$ be the $h$-vector of a level algebra $A$ of type 2 and socle degree $e$. Then, for every Gorenstein quotient $B$ of $A$ having socle degree $e$, the reverse of $h$ can be written as the sum of two $O$-sequences (which are $h$-vectors of quotients of $A$), the longer of these $O$-sequences being (the reverse of) the $h$-vector of $B$.\\}
\\\indent 
{\bf Proof.} See $[GHMS]$, Corollary 2.11.{\ }{\ }\qed \\
\\\indent 
The following is a stronger version of $[GHMS]$'s decomposition result:\\
\\\indent
{\bf Theorem 2.3.} {\it Let $h$ be the $h$-vector of a level algebra $A$ of type 2, and let $M$ be the Inverse System module associated to $A$. Then, for every choice of two forms $F$ and $G$ generating $M$, $h$ can be written as $$h=h^{'}+h^{''}+h^{'''},$$ where $h^{''}$ is an $O$-sequence, the reverses of $h^{'}$ and $h^{'''}$ are $O$-sequences, and $h^{'}+h^{''}$ and $h^{''}+h^{'''}$ are the Gorenstein $h$-vectors given by $\langle F \rangle $ and $\langle G \rangle $. The three $O$-sequences above are $h$-vectors of quotients of $A$.\\}
\\\indent 
{\bf Proof.} Let $M=\langle F,G\rangle \subseteq S$ be an Inverse System module giving the $h$-vector $h$ (i.e. $h$ is the $h$-vector of $A=R/M^{-1}$). Notice that $A=R/(\langle F \rangle^{-1}\cap \langle G\rangle^{-1})$. By standard facts about Inverse Systems, $N=\langle F\rangle \cap \langle G\rangle $ is the finitely generated submodule of $S$ corresponding to the algebra $R/(\langle F\rangle^{-1}+\langle G\rangle^{-1})$ (which is clearly a quotient of $A$), say with $h$-vector $h^{''}$. By the exclusion-inclusion principle it is easy to see that, for each $d\leq e$, $$\dim_k (\langle F\rangle_d/N_d)+\dim_k (\langle G\rangle_d/N_d)+ \dim_kN_d=\dim_k M_d.$$ \indent Furthermore, $\langle F\rangle_d/N_d\cong M_d/\langle G\rangle_d$ as $k$-vector spaces, and, similarly, $\langle G\rangle_d/N_d\cong M_d/\langle F \rangle_d$. By Theorem 2.1, the reverse of the sequence $h^{'}$ of the dimensions of the $M_d/\langle G\rangle_d$ is an $O$-sequence (and is also the $h$-vector of a quotient of $A$), and the same is true for the reverse of $h^{'''}$, the sequence of the dimensions of the $M_d/\langle F\rangle_d$. Clearly, by definition, $h^{'}+h^{''}$ is the Gorenstein $h$-vector given by $\langle F\rangle $, and $h^{''}+h^{'''}$ is the Gorenstein $h$-vector given by $\langle G \rangle $. This concludes the proof of the theorem.{\ }{\ }\qed \\
\\\indent 
{\bf Remark 2.4.} i). We have stated Theorem 2.3 in the most interesting case, that is for $h_e=2$, but the proof above can be extended naturally to any $h_e$: in fact it can be shown in a similar fashion that the $h$-vector of a level algebra $A=R/M^{-1}$ of type $h_e$ (for any choice of $h_e$ forms $F_1,...,F_{h_e}$ generating $M$) can be written as the sum of $2^{h_e}-1$ suitable vectors, such that $2^{h_e}-1-h_e$ of them are $O$-sequences and the reverses of the remaining $h_e$ are also $O$-sequences.\\\indent
ii). There are cases in which Theorem 2.3 applies but Corollary 2.2 does not. Consider, for instance, $h=(1,3,6,10,9,7,5,2)$. We want to show that $h$ is not a level $h$-vector. (Notice that the reverse of $h$ can be written as $(1,3,4,5,5,4,3,1)+(0,0,2,5,4,3,2,1)$, which are both $O$-sequences, and therefore Corollary 2.2 does not suffice to rule out $h$.)\\\indent
Suppose that $A=R/M^{-1}$ is a level algebra having $h$-vector $h$. Consider any form $H\in M$ of degree 7 having three first derivatives, and let $h^{H}$ be the $h$-vector given by $\langle H \rangle $. It is easy to check (see Stanley's Proposition 2.7 below) that, in order to let the reverse of $h-h^{H}$ be an $O$-sequence, the only possibility for $h^{H}$ is $h^{H}=(1,3,4,5,5,4,3,1)$. Moreover, since the entry of degree 3 of $h$ is equal to 10, there cannot be a form of degree 7 in $M$ having only 2 first derivatives. In fact, the $h$-vector given by such a form should be clearly bounded from above by $(1,2,3,4,4,3,2,1)$, but $5+4<10$. Thus, every form of degree 7 in $M$ must give the $h$-vector $(1,3,4,5,5,4,3,1)$. Therefore, it easily follows by Theorem 2.3 that, for any choice of two forms $F$ and $G$ generating $M$, $h$ must decompose as $$h=h^{'}+h^{''}+h^{'''}=(0,0,2,5,4,3,2,1)+(1,3,2,0,1,1,1,0)+(0,0,2,5,4,3,2,1),$$ but this is a contradiction, since $h^{''}$ is not an $O$-sequence. Hence $h$ is not a level $h$-vector.\\
\\\indent 
The following result of Iarrobino is a fundamental tool we will use in this paper. In particular, given a level algebra $A$ of type 2, Theorem 2.5 guarantees the existence of a Gorenstein quotient of $A$ of the same socle degree having at least a certain codimension.\\
\\\indent 
{\bf Theorem 2.5} (Iarrobino). {\it Let $h=(1,h_1,...,h_{e-1},h_e=2)$ be the $h$-vector of a level algebra $A$ of type 2, and fix an integer $u$ such that $1\leq u\leq e$. Suppose that, for some integer $\delta_{u}\geq 0$, $h_{e-u}\geq 2h_u-2-3\delta_{u}$. Then there exists a Gorenstein quotient $B$ of $A$ of socle degree $e$, having $h$-vector $h^{'}=(1,h^{'}_1,...,h^{'}_e=1)$ and such that $h^{'}_{u}\geq h_u-\delta_{u}$.\\}
\\\indent 
{\bf Proof.}  See $[Ia]$, Theorem 2.4 (where a stronger conclusion is stated).{\ }{\ }\qed \\
\\\indent
{\bf Remark 2.6.} i). (Unfortunately) Iarrobino's lower-bound is a very good one! Let us consider a level $h$-vector of the form $(1,r,...,r,2)$, and let $u=1$. The least integer $\delta_1$ satisfying the inequality $r\geq 2r-2-3\delta_1$ is clearly $\delta_1 =\lceil {r-2\over 3}\rceil ,$ where, as usual, $\lceil \alpha \rceil $ indicates the least integer greater than or equal to $\alpha $.\\\indent
For every $r\geq 2$, we will exhibit a level algebra $A$ of socle degree $e$, having an $h$-vector $h=(1,h_1,...,h_e)$ of the form $(1,r,...,r,2)$ and such that every Gorenstein quotient of $A$ of socle degree $e$ has codimension $h^{'}_1\leq r-\lceil {r-2\over 3}\rceil $. As a consequence we have that, in this case, Iarrobino's lower-bound is sharp.\\\indent
Let $r$ be a multiple of 3, say $r=3p$ for some positive integer $p$. Consider the Inverse System module $M=\langle F,G \rangle $, where $F=y_{p+1}y_1^{e-1}+y_{p+2}y_2^{e-1}+...+y_{2p}y_p^{e-1}$ and $G=y_{2p+1}y_1^{e-1}+y_{2p+2}y_2^{e-1}+...+y_{3p}y_p^{e-1}$. It is easy to see that the algebra $A$ associated to $M$ satisfies our requirements: in fact, the $h$-vector of $A$ is clearly of the form $(1,r,...,r,2)$ (indeed, it is equal to $(1,r,r,...,r,r,2)$), and every Gorenstein quotient of $A$ of socle degree $e$ has codimension $h^{'}_1=2p=r-\lceil {r-2\over 3}\rceil $. This last fact is true since the Gorenstein quotients of $A$ of socle degree $e$ correspond to the Inverse System modules generated by $\mu F+\lambda G=(\mu y_{p+1}+\lambda y_{2p+1})y_1^{e-1}+...+(\mu y_{2p}+\lambda y_{3p})y_p^{e-1}$, where $[\mu :\lambda ]$ ranges in ${\bf P}^1(k)$, and, by changing variables, it is immediate to see that all of these forms have exactly $2p$ first derivatives.\\\indent
The cases $r=3p+1$ and $r=3p+2$ can be treated similarly. If $r=3p+1$, consider $M=\langle F,G\rangle $, where $F=y_{p+1}y_1^{e-1}+y_{p+2}y_2^{e-1}+...+y_{2p}y_p^{e-1}+y_{3p+1}^e$ and $G=y_{2p+1}y_1^{e-1}+y_{2p+2}y_2^{e-1}+...+y_{3p}y_p^{e-1}$. Then the $h$-vector of $A=R/M^{-1}$ has the form $(1,r,...,r,2)$ (again, it is equal to $(1,r,r,...,r,r,2)$), and every Gorenstein quotient of $A$ of socle degree $e$ has codimension $h^{'}_1=2p+1=r-\lceil {r-2\over 3}\rceil $, except for $R/\langle G\rangle^{-1}$, that has codimension $2p$.\\\indent
If $r=3p+2$, consider $M=\langle F,G \rangle $, where $F=y_{p+1}y_1^{e-1}+y_{p+2}y_2^{e-1}+...+y_{2p}y_p^{e-1}+y_{3p+1}^e$ and $G=y_{2p+1}y_1^{e-1}+y_{2p+2}y_2^{e-1}+...+y_{3p}y_p^{e-1}+y_{3p+2}^e$. Then the $h$-vector of $A=R/M^{-1}$ has the form $(1,r,...,r,2)$ (as in the previous two cases, it is equal to $(1,r,r,...,r,r,2)$), and every Gorenstein quotient of $A$ of socle degree $e$ has codimension $h^{'}_1=2p+2=r-\lceil {r-2\over 3}\rceil $, except for $R/\langle F\rangle^{-1}$ and $R/\langle G\rangle^{-1}$, that have both codimension $2p+1$.\\\indent
ii). Actually, the examples above show that Iarrobino's lower-bound is sharp for every $u=1,...,e$. The proof for the other values of $u$ is similar to that we have given for $u=1$. We just notice that, in the example above, when $r$ is a multiple of 3 there is exactly one possible $h$-vector for the Gorenstein quotients of $A$ of socle degree $e$, and this, for every $u$, $2\leq u\leq e-1$, is equal, in the $u$-th place, to $r-\lceil {r-2\over 3}\rceil $, where $\lceil {r-2\over 3}\rceil $ is the least integer $\delta_u$ satisfying Iarrobino's inequality $h_{e-u}\geq 2h_u-2-3\delta_{u}$ (i.e. $r\geq 2r-2-3\delta_u$).\\\indent
In each of the other two cases (i.e. when $r=3p+1$ and $r=3p+2$), there are exactly two possible $h$-vectors for the Gorenstein quotients of $A$ of socle degree $e$ (one less than the other in every place), and the larger of these, for every $u$, $2\leq u\leq e-1$, is equal, in the $u$-th place, to $r-\lceil {r-2\over 3}\rceil $, where $\lceil {r-2\over 3}\rceil $ is again the least integer $\delta_u$ satisfying Iarrobino's inequality.\\\indent
iii). A remarkable fact in cases like those we described above is that there is a considerable difference between the codimension $r$ of the type two level algebra and the codimension $h^{'}_1$ of any of its Gorenstein quotients having the same socle degree, since $h^{'}_1\leq r-\lceil {r-2\over 3}\rceil ={2\over 3}r +O(1)$.\\
\\
\indent Recall that a vector $v=(1,v_1,v_2,...,v_d)$ is said to be {\it differentiable} if its first difference, $$\Delta v=((\Delta v)_0=1,(\Delta v)_1=v_1-1,(\Delta v)_2=v_2-v_1,...,(\Delta v)_d=v_d-v_{d-1}),$$ is an $O$-sequence.\\\indent
An $h$-vector $h=(1,h_1,...,h_e)$ is an {\it SI-sequence} if it is symmetric with respect to ${e\over 2}$ and if its first half, $(1,h_1,...,h_{\left\lfloor {e\over 2}\right\rfloor })$, is differentiable (as usual, $\lfloor \alpha \rfloor $ denotes the greatest integer less than or equal to $\alpha $).\\
\\\indent 
{\bf Proposition 2.7} (Stanley). {\it Let $h=(1,h_1,...,h_e)$, $h_1\leq 3$. Then $h$ is a Gorenstein $h$-vector if and only if $h$ is an SI-sequence.\\}
\\\indent 
{\bf Proof.} See $[St]$, Theorem 4.2. For $h_1=2$ this result was already known to Macaulay (cf. $[Ma]$).{\ }{\ }\qed \\
\\\indent 
The following result determines the Gorenstein $h$-vectors of codimension 4 having the second entry $h_2\leq 7$. (For $8\leq h_2\leq 10$ the problem of characterizing these $h$-vectors remains a very interesting open one.)\\
\\\indent 
{\bf Proposition 2.8} (Iarrobino-Srinivasan). {\it Let $h=(1,h_1,h_2,...,h_e)$, with $h_1=4$ and $h_2\leq 7$. Then $h$ is a Gorenstein $h$-vector if and only if $h$ is an SI-sequence.\\}
\\\indent 
{\bf Proof.} See $[IS]$, Theorem 3.2, Corollary 3.3 and Proposition 3.4.{\ }{\ }\qed \\
\\\indent 
{\bf Theorem 2.9.} {\it Let $r\leq 4$. Then $h=(1,r,...,r,2)$ is a level $h$-vector if and only if $h$ can be written as $$h=h^{(r-1)}+(0,1,1,...,1),$$ where $h^{(r-1)}$ is a Gorenstein $h$-vector of codimension $r-1$.\\}
\\\indent 
{\bf Proof.} \lq \lq $\Longleftarrow $": Let $e$ be the socle degree of $h$ and let $F=F(y_1,...,y_{r-1})\in S=k[y_1,...,y_r]$ be a form of degree $e$ which generates the Inverse System module giving $h^{(r-1)}$. Then $h$ is clearly the $h$-vector given by the module $\langle F,y_r^e\rangle $.\\
\indent \lq \lq $\Longrightarrow $": The case $h=(1,r,r,2)$ is trivial. Thus, suppose for the rest of the proof that the socle degree $e$ of $h$ is at least 4. Let $A$ be a level algebra having $h$-vector $h$. Since $r\leq 4$, by Theorem 2.5 it is easy to see (since for $u=1$ we can choose $\delta_{u}=1$) that there exists a Gorenstein quotient of $A$ having socle degree $e$ and codimension at least $r-1$ (for $r=2$ we can always obtain codimension $r$).\\
\indent Let $3\leq r\leq 4$, and suppose that there are no Gorenstein quotients of $A$ having socle degree $e$ and codimension $r$. Hence there exists a Gorenstein quotient $B$ of socle degree $e$ and codimension $r-1$. By Corollary 2.2, the reverse of the difference between the $h$-vector of $A$ and that of $B$ is an $O$-sequence. Since this $O$-sequence must begin with $(1,1,...)$ and end with $(...,1,0)$, by Macaulay's Theorem 1.2, it has to have the form $(1,1,...,1,1,0)$, and the theorem follows.\\\indent 
Suppose then that there exists a Gorenstein quotient of $A$ having socle degree $e$ and codimension $r$, $2\leq r\leq 4$. Let $h=(1,r,...,b,r,2)$ be the $h$-vector given by the Inverse System module $M$ associated to $A$. Hence there is a form $F\in M$ of degree $e$ having $r$ linearly independent first derivatives. Clearly, by Corollary 2.2, $h$ is equal to the Gorenstein $h$-vector of $\langle F \rangle $ plus $(0,0,...,0,1)$.\\\indent 
Now, let us assume for the moment that $b\leq {r\choose 2}+1$ (this fact will be shown in Theorem 3.5, i)).\\\indent 
If $r=2$, then we are done: in fact, by Macaulay's Theorem 1.2, we cannot have $b<2$, so $b=2$, whence it clearly follows that $h=(1,2,2,...,2,2,2)=(1,1,1,...,1,1,1)+(0,1,1,...,1,1,1)$.\\\indent 
Let $r=3$; thus $b\leq 4$. Hence, the $h$-vector of the Gorenstein quotient of $A$ corresponding to $\langle F \rangle $, by Proposition 2.7, is equal to $h^{'}=(1,3,3,...,3,3,1)$ or $h^{''}=(1,3,4,...,4,3,1)$, where the first difference of the first half of $h^{''}$ has the form $(1,2,1,1,...,1,0,...,0)$. This clearly implies that the $h$-vector given by $\langle F \rangle $ is the sum of a Gorenstein $h$-vector of codimension 2 and $(0,1,1,...,1,0)$, and the result easily follows.\\\indent 
Let $r=4$. We have $b\leq 7$. As in the case $r=3$, it will be enough to show that a Gorenstein $h$-vector of the form $h^{'}=(1,4,b,...,4,1)$, with $b\leq 7$, is the sum of a Gorenstein $h$-vector of codimension 3 and $(0,1,1,...,1,0)$. By Proposition 2.8, the first difference of the first half of $h^{'}$, that is $\Delta^{'}=(1,3,b-4,...)$, is an $O$-sequence. But $\Delta^{'}$, from degree 2 on, is clearly equal to the first difference of the first half of $h^{'}-(0,1,1,...,1,0)$, that is to $\Delta^{''}=(1,2,b-4,...)$. Since $b-4\leq 3$, it follows from Macaulay's theorem that $\Delta^{''}$ is also an $O$-sequence. Hence, by Proposition 2.7, $h^{'}-(0,1,1,...,1,0)$ is a Gorenstein $h$-vector of codimension 3. This concludes the proof of the theorem.{\ }{\ }\qed\\
\\\indent 
{\bf Remark 2.10.} i). Notice that the proof above cannot be extended to $r\geq 5$, since we do not have enough information about the Gorenstein $h$-vectors of codimension larger than 3. Therefore the characterization of level $h$-vectors of the form $(1,r,...,r,2)$ for $r\geq 5$ remains an open problem.\\\indent
ii). Theorem 2.9 shows, in particular, that, for $r\leq 4$, level $h$-vectors of the form $(1,r,...,r,2)$ are unimodal with their first half being differentiable, providing a positive answer, for these $h$-vectors, to a question of Geramita {\it et al.} ($[GHMS]$, Question 4.4).\\
\\
\section{Upper-bounds for the entries of a level $h$-vector of type 2}
\indent 

In this section we investigate the intervals in which the entries of a level $h$-vector of type 2 may range. The main result is (under certain conditions) an (entry by entry) sharp upper-bound for the level $h$-vectors of type 2, given the next to last entry.\\
\\\indent
{\bf Lemma 3.1} (Iarrobino). {\it Let $F=\sum_{t=1}^m L_t^d$ be a form of degree $d$ in $S=k[y_1,...,y_r]$, where the $L_t=\sum_{k=1}^rb_{tk}y_k$ are linear forms, and let $I_{\bf b}=Ann(F)\subseteq R$, with ${\bf b}=(b_{11},...,b_{mr})$. Then there exists a non-empty open subseteq $U$ of $k^{mr}$ such that, for every ${\bf b}\in U$, the Gorenstein artinian algebras $R/I_{\bf b}$ all have the same $h$-vector, denoted by: $$h(m,d)=(1,h_1(m,d),...,h_d(m,d)=1),$$ where, for $j=1,...,d$, $$h_j(m,d)=\min \lbrace m,\dim_kR_j,\dim_kR_{d-j}\rbrace .$$}
\\
\indent {\bf Proof.} See $[Ia2]$, Proposition 4.7.{\ }{\ }\qed \\
\\\indent 
{\bf Lemma 3.2} (Iarrobino). {\it Let $h=(1,h_1,...,h_d)$ be the $h$-vector of a level algebra $A=R/I$, where $I$ annihilates the $R$-submodule $M$ of $S$. Let $m\leq {r-1+d\choose d}-h_d$. Then, if $F$ is the sum of the $d$-th powers of $m$ generic linear forms (the Gorenstein $h$-vector of $R/\langle F\rangle^{-1}$ is $h(m,d)$, given by Lemma 3.1), then the level algebra associated to $M^{'}=\langle M,F \rangle $ has $h$-vector $H=(1,H_1,...,H_d),$ where, for $i=1,...,d$, $$H_i=\min \lbrace h_i+h_i(m,d),{r-1+i\choose i}\rbrace .$$}
\\\indent 
{\bf Proof.} See $[Ia2]$, Theorem 4.8 A, where this result appears under more general hypotheses.{\ }{\ }\qed \\
\\\indent 
Let us consider the next to last entry $a$ of a level $h$-vector $h$ of type 2 and codimension $r$. Using Iarrobino's Theorem 2.5, we are able to supply a sharp lower-bound ($a=r$) for $r\leq 7$. In a certain sense, this result cannot be improved: in fact, Brian Coolen, a doctoral student of Professor A.V. Geramita (personal communication), has produced some examples that show that, for $r>7$, $a=r$ is no longer a lower-bound. It would be interesting to determine a sharp lower-bound for $a$ for any codimension $r$.\\
\\\indent 
{\bf Theorem 3.3.} {\it Let $h=(1,r,...,a,2)$ be a level $h$-vector of type 2. Then:\\
\indent i). $a\leq 2r$, and $a$ may assume any integral value in the range $[r,2r]$.\\
\indent ii). If $r\leq 7$, then $a\geq r$.\\}
\\\indent 
{\bf Proof.} i). Let $M=\langle F,G \rangle $ be the Inverse System module associated to a level algebra $A$ having $h$-vector $h$, where $\deg (F)=\deg (G)=e$. Obviously, $a\leq 2r$, since $F$ and $G$ can have at most $r$ first derivatives each.\\\indent
In order to obtain any integral value of $a$ in $[r,2r]$, write $a=d_1+d_2$, with $1\leq d_1,d_2\leq r$. By Lemmata 3.2 and 3.1, it clearly suffices to choose $G$ as the sum of powers of $d_1$ generic linear forms and $F$ as the sum of powers of $d_2$ generic linear forms.\\
\indent ii). We will make use of Iarrobino's Theorem 2.5 with $u=1$. Here $h_1=r$ and $h_{e-1}=a$. Suppose $a<r$. Let us first show that the codimension $h_1^{'}$ of any Gorenstein quotient $B$ of $A$ of socle degree $e$ must satisfy the inequality $h_1^{'}<a-1$. In fact, $h_1^{'}>a$ is obviously impossible. Furthermore, if $h_1^{'}\leq a$, then the reverse of the difference between $(1,r,...,a,2)$ and $(1,h_1^{'},...,h_1^{'},1)$ has the form $(1,a-h_1^{'},...,r-h_1^{'},0)$, which must be an $O$-sequence by Corollary 2.2. But it is immediate to see that, for $a<r$, this is impossible if $h_1^{'}\in \lbrace a-1,a\rbrace $. Hence $h_1^{'}<a-1$, as we desired.\\
\indent By Theorem 2.5, if $a\geq 2r-2-3\delta $ for some integer $\delta \geq 0$, then, for some Gorenstein quotient $B$ of $A$ having codimension $h_1^{'}$, we have $h_1^{'}\geq r-\delta $. Hence, in order to show that some value of $a$ satisfies the inequality $h_1^{'}\geq a-1$, it suffices to show that $r-\delta \geq a-1$ or, by Theorem 2.5, that \begin{equation}\label{aa} r-\lceil {2r-2-a\over 3}\rceil -a+1\geq 0.\end{equation}
\indent It is easy to check that, if $a$ satisfies (\ref{aa}), then also $a-1$ does. Therefore, in order to show that $a$ cannot be less than $r$, it is enough to show that (\ref{aa}) holds for $a=r-1$. But, for $a=r-1$, (\ref{aa}) becomes $$r-\lceil {2r-2-(r-1)\over 3}\rceil-(r-1)+1\geq 0,$$ i.e. $\lceil {r-1\over 3}\rceil \leq 2$, which is true if (and only if) $r\leq 7$. Hence, for $r\leq 7$, we have $a\geq r$, as we wanted to show.{\ }{\ }\qed \\
\\\indent
The following lemma is a combinatorial result that we will need in the next two theorems in order to determine an upper-bound for the entries of our level $h$-vectors of type 2.\\
\\\indent 
{\bf Lemma 3.4.} {\it Let $r$, $i$ and $a$ be integers such that $i\geq 1$ and $2\leq r\leq a\leq 2r$. Then: $$\max \lbrace {r_1-1+i\choose i}+{r_2-1+i\choose i} {\ } \mid {\ } 1\leq r_1,r_2\leq r,{\ } r_1+r_2=a\rbrace $$$$
=\cases{{r-2+i\choose i}+1,&if $a=r$\cr
	{r-1+i\choose i}+{a-r-1+i\choose i},&if $r<a\leq 2r$.\cr}
$$}
\\\indent 
{\bf Proof.} Notice that ${r_1-1+i\choose i}+{r_2-1+i\choose i}$, for $r_1$ and $r_2$ as in the hypotheses, takes on the value of the statement when the absolute value of $r_1-r_2$ is as large as possible. So, it remains to show that the maximal value of the sum above occurs when the absolute value of $r_1-r_2$ is as large as possible.\\\indent
Now, since $r_1-r_2=(r_1-1+i)-(r_2-1+i)$ and $r_1+r_2$ is always fixed to equal $a$, it suffices to show that, for $r_1\geq r_2\geq 2$, \begin{equation}\label{com}{r_1-1+i\choose i}+{r_2-1+i\choose i}\leq {r_1+i\choose i}+{r_2-2+i\choose i},$$ i.e. that $${r_2-1+i\choose i}-{r_2-2+i\choose i}\leq {r_1+i\choose i}-{r_1-1+i\choose i}.\end{equation}\indent
By the standard Pascal triangle inequality, (\ref{com}) becomes $${r_2-2+i\choose i-1}\leq {r_1-1+i\choose i-1},$$ which is true if and only if $r_2-2+i\leq r_1-1+i$, and this is the case since we are assuming $r_2\leq r_1$. This proves the lemma.{\ }{\ }\qed\\
\\\indent 
The next theorem gives information on the third-last entry $b$ of a level $h$-vector $h=(1,r,...,b,a,2)$ of type 2 and codimension $r$, in particular supplying a sharp upper-bound for $b$.\\
\\\indent 
{\bf Theorem 3.5.} {\it Let $h=(1,r,...,b,a,2)$ be a level $h$-vector of socle degree $e$. Then:\\\indent 
i). If $a=r$, then $b\leq {r\choose 2}+1$. Moreover, $(...,b,r,2)$ is the final part of a level $h$-vector of codimension $r$ for any integer $b\in [r,{r\choose 2}+1]$.\\\indent 
ii). If $r<a\leq 2r$, then $b\leq \min \lbrace {r+1\choose 2}+{a-r+1\choose 2},{r+e-3\choose e-2}\rbrace $. Moreover, $(...,b,a,2)$ is the final part of a level $h$-vector of codimension $r$ for any integer $b\in [a,\min \lbrace {r+1\choose 2}+{a-r+1\choose 2},{r+e-3\choose e-2}\rbrace ]$.\\}
\\\indent 
{\bf Proof.} i). Let $I^{-1}=M=\langle F,G\rangle $ be the Inverse System module associated to a level algebra $A=R/I$ having $h$-vector $h$, where $\deg (F)=\deg (G)=e$. Suppose that there exists a form $H$ of degree $e$ in $M$ having $r_H<r$ linearly independent first derivatives. Then the $h$-vector given by $\langle H\rangle $ starts with $(1,r_H,b_H,...)$, with $b_H\leq {r_H+1\choose 2}$ by Macaulay's Theorem 1.2. By Corollary 2.2, $(2,r,b,...)-(1,r_H,b_H,...)=(1,r-r_H,b-b_H,...)$ is an $O$-sequence, and therefore $b-b_H\leq {r-r_H+1\choose 2}$. Thus, $b\leq {r_H+1\choose 2}+{r-r_H+1\choose 2}$. Hence, by Lemma 3.4 for $i=2$ and $a=r$, we have $b\leq {r\choose 2}+1,$ as we wanted to show.\\\indent 
Therefore, from now on, let us suppose that all the forms of degree $e$ in $M$ have $r$ linearly independent first derivatives. Since $a=r$, the first derivatives of $G$ are generated by those of $F$, and therefore, for $i=1,...,r$, we can write \begin{equation}\label{g}G_{y_i}=c_{1,i}F_{y_1}+...+c_{r,i}F_{y_r}, {\ }{\ }{\ }c_{j,i}\in k,\end{equation} where, for convenience, we write $F_{y_i}$ for the partial derivative of the form $F$ with respect to $y_i$.\\\indent 
Notice that the matrix $C =(c_{i,j})_{i,j=1,...,r}$ is invertible, since there is no dependence relation among the first derivatives of $G$, i.e. $\det{C}\neq 0$.\\\indent
Suppose, for the moment, that $C$ is not diagonal. Then, without loss of generality, assume that the first row is different from $(c_{1,1},0,...,0)$, i.e. $c_{1,t}\neq 0$ for some $2\leq t\leq r$. Since $G_{y_1,y_i}=G_{y_i,y_1}$, we have, for each $i$, $2\leq i\leq r$, \begin{equation}\label{1i}c_{1,1}F_{y_1,y_i}+c_{1,2}F_{y_2,y_i}+...+c_{1,r}F_{y_r,y_i}-c_{i,1}F_{y_1,y_1}-c_{i,2}F_{y_2,y_1}-...-c_{i,r}F_{y_r,y_1}=0.\end{equation}
\indent It is easy to see that these $r-1$ equations (\ref{1i}) are independent (as linear relations among the second derivatives of $F$), since, for any $2\leq j\leq r$, the term $F_{y_t,y_j}$, whose coefficient is $c_{1,t}\neq 0$, is only in one of the equations (\ref{1i}), namely the one for which $i=j$.\\
\indent Thus, $b\leq {r+1\choose 2}-(r-1)={r\choose 2}+1$, as we desired.\\\indent 
Now suppose $C=$diag$\lbrace c_1,...,c_r\rbrace $. We have $G_{y_i}=c_iF_{y_i}$. Hence the equality $G_{y_i,y_j}=G_{y_j,y_i}$ implies, for $1\leq i,j\leq r$, that $(c_i-c_j)F_{y_i,y_j}=0$. If $c_1\neq c_i$ for every $i\geq 2$ then $F_{y_1,y_i}=0$ for $i\geq 2$, whence we have again the desired upper-bound $b\leq {r+1\choose 2}-(r-1)={r\choose 2}+1$.\\\indent 
Hence suppose $c_1=c_2$. If $r>2$, let $c_1=c_2$ be different from $c_i$ for every $i\geq 3$; thus we have $2(r-2)$ second derivatives of $F$ equal to 0. Since $2(r-2)\geq r-1$, again we obtain the desired upper-bound for $b$. By induction, suppose that, for some $i<r$, $c_1=c_2=...=c_i=c$, and that $c\neq c_j$ for every $j$, $i+1\leq j\leq r$. Hence $i(r-i)$ second derivatives of $F$ are equal to 0 and, since $i(r-i)\geq r-1$ for $i<r$, the upper-bound for $b$ is again satisfied.\\\indent 
Therefore let $c_1=...=c_r=c$. Hence $G_{y_i}=cF_{y_i}$ for $i=1,...,r$, i.e. all the first derivatives of $G-cF$ are equal to 0. Thus the form $G-cF$ is constant, and therefore is equal to 0, but this is a contradiction, since $F$ and $G$ were supposed to be linearly independent.\\\indent 
This completes the proof that, for $a=r$, $b\leq {r\choose 2}+1$.\\\indent 
In order to show that $(...,b,r,2)$ is the final part of a level $h$-vector of codimension $r$ for every $b\in [r,{r\choose 2}+1]$, it is enough to consider the Inverse System module $M^{'}=\langle F^{'},G^{'}\rangle $, where we choose $F^{'}$ only in variables $y_1,...,y_{r-1}$ to be the sum of the $e$-th powers of $b-1$ generic linear forms and $G^{'}=y_r^e$. Then, by Lemma 3.1, $M^{'}$ gives the desired $h$-vector.\\\indent 
ii). Let $h=(1,r,...,b,a,2)$ be as in the statement, and let $A=R/I$ be a level algebra having $h$-vector $h$, with $I^{-1}=M=\langle F,G\rangle $. The $h$-vector given by $\langle F\rangle $ ends with $(...,b_F,r_F,1)$, where $b_F\leq {r_F+1\choose 2}$. By Corollary 2.2, $(2,a,b,...)-(1,r_F,b_F,...)=(1,a-r_F,b-b_F,...)$ is an $O$-sequence, and therefore $b-b_F\leq {a-r_F+1\choose 2}$. Thus $b\leq {r_F+1\choose 2}+{a-r_F+1\choose 2}$. Hence, by Lemma 3.4 for $i=2$, we have $b\leq {r+1\choose 2}+{a-r+1\choose 2}.$\indent 
Of course $b\leq {r+e-3\choose e-2}$, since this is the dimension of $R_{e-2}$ as a $k$-vector space. This proves the upper-bound for $b$.\\\indent 
In order to show that $(...,b,a,2)$ is the final part of a level $h$-vector of codimension $r$ for every $b\in [a,\min \lbrace {r+1\choose 2}+{a-r+1\choose 2},{r+e-3\choose e-2}\rbrace ]$, consider the Inverse System module $M^{'}=\langle F^{'},G^{'}\rangle $, where we choose $G^{'}$ only in variables $y_1,...,y_{a-r}$ and $F^{'}$ in all of the $r$ variables, both of the same degree $e$. Let us make use of Lemmata 3.2 and 3.1. For $b=a$ choose $G^{'}$ to be the sum of the $e$-th powers of $a-r$ generic linear forms and $F^{'}$ the sum of the $e$-th powers of $r$ generic linear forms. For $b=a+1$ do the same with $r+1$ linear forms (instead of $r$) for $F^{'}$ and the same number, $a-r$, for $G^{'}$. We continue this procedure up to $b=\min \lbrace {r+1\choose 2}+(a-r),{r+e-3\choose e-2}\rbrace $ (clearly, as soon as the minimum is ${r+e-3\choose e-2}$ the whole process finishes). Instead, from $b={r+1\choose 2}+(a-r+1)$ up to ${r+e-3\choose e-2}$, we now increase the number of generic linear forms for $G^{'}$ by 1 at the time, letting $F^{'}$ remain always the sum of ${r+1\choose 2}$ $e$-th powers. In this way, by Lemmata 3.2 and 3.1, for the third-last entry $b$ of the $h$-vector $(...,b,a,2)$ given by $M^{'}$ we clearly obtain every value between $a$ and $\min \lbrace {r+1\choose 2}+{a-r+1\choose 2},{r+e-3\choose e-2}\rbrace$, and the proof of the theorem is complete.{\ }{\ }\qed \\
\\\indent 
{\bf Remark 3.6.} The upper-bound we have shown for $b$ when $a=r$ was already known to S\"oderberg ($[So]$), who obtained the same result using a different approach.\\
\\\indent
{\bf Example 3.7.} Theorem 3.5 does not give information for $b<a$, even if there exist such cases. For instance, consider the Gorenstein $h$-vector $h=(1,52,51,52,51,52,1)$ (its existence is shown in $[Bo]$, Example 2.11). If we add the 6-th power of a generic linear form to an Inverse System module giving this $h$-vector, by Lemma 3.2 we obtain the level $h$-vector $(1,52,52,53,52,53,2)$. Theorem 3.5, instead, guarantees the existence of level $h$-vectors of the form $(1,52,...,b,53,2)$ for every integer $b\in [53,\min \lbrace {52+1\choose 2}+1,{55\choose 4}\rbrace ]=[53,1378]$, where $1378$ is also an upper-bound.\\
\\\indent 
The next result is generalization of Theorem 3.5: given the last two entries $(a,2)$, $r\leq a\leq 2r$, of a level $h$-vector $h$ of type 2 and codimension $r$, we provide a sharp upper-bound for all of $h$ (we need to assume that $r\leq 5$ if $a=r$).\\
\\\indent 
{\bf Theorem 3.8.} {\it Let $h=(1,r,...,h_{e-i},...,a,2)$ be a level $h$-vector of socle degree $e$, where we fix the index $i$, $2\leq i\leq e-2$. Then:\\
\indent i). Let $2\leq r\leq 5$. If $a=r$, then $$h_{e-i}=h_i\leq \min \lbrace {r-2+i\choose i}+1,{r-2+e-i\choose e-i}+1\rbrace .$$ \indent Moreover, there exists a level $h$-vector of the form $(1,r,...,h_{e-i},...,r,2)$ for any integer $h_{e-i}\in [r,\min \lbrace {r-2+i\choose i}+1,{r-2+e-i\choose e-i}+1\rbrace ].$\\
\indent ii). Let $r\geq 2$. If $r<a\leq 2r$, then $$h_{e-i}\leq \min \lbrace {r-1+i\choose i}+{a-r-1+i\choose i},{r-1+e-i\choose e-i}\rbrace .$$ \indent Moreover, there exists a level $h$-vector of the form $(1,r,...,h_{e-i},...,a,2)$ for any integer $h_{e-i}\in [a, \min \lbrace {r-1+i\choose i}+{a-r-1+i\choose i},{r-1+e-i\choose e-i}\rbrace ].$\\
\indent iii). For any $a\in[r,2r]$ (always assuming $r\leq 5$ if $a=r$), the upper-bounds for the $h_{e-i}$'s given above can all be assumed simultaneously. In other words, the $h$-vector $H=(1,r,h_2,...,h_{e-2},a,2)$ given by the maximum value of the $h_{e-i}$'s for each $i=2,...,e-2$ is the (entry by entry) maximum for all the level $h$-vectors of the form $(1,r,...,a,2)$.\\}
\\\indent 
{\bf Proof.} i). Let $a=r$. By Inverse Systems, let us suppose that $h$ be given by an $R$-submodule $M=\langle F,G \rangle $ of $S$, where $\deg{F}=\deg{G}=e$, and that $A$ is the algebra associated to $M$. By Iarrobino's Theorem 2.5, since $r\leq 5$, it is easy to see that there exists a Gorenstein quotient $B$ of $A$ having codimension at least $r-1$. Hence, by Corollary 2.2, the reverse of $h$ can be written as the sum of the Gorenstein $h$-vector of $B$ and $(1,1,...,1,0)$ (if $B$ has codimension $r-1$) or $(1,0,0,...,0)$ (if $B$ has codimension $r$). In either case, clearly, $h_i=h_{e-i}$ for $i=1,...,e-1$.\\\indent 
Suppose that there exists no Gorenstein quotient of $A$ of codimension $r$. Hence there exists a form $F\in M$ of degree $e$ giving a Gorenstein $h$-vector of codimension $r-1$. Therefore, the $h$-vector of $\langle F\rangle $ starts with $(1,r-1,...)$, and thus its entry of degree $i$, $1\leq i\leq e-1$, by Macaulay's Theorem 1.2 is at most ${r-2+i\choose i}$. By Corollary 2.2 and the fact that $h_i=h_{e-i}$, the upper-bound for $h_{e-i}$ easily follows.\\\indent 
Suppose now that there exists a form $F\in M$ of degree $e$ having $r$ linearly independent first derivatives. Since, by Theorem 3.5, i), the Gorenstein $h$-vector given by $F$ starts with $(1,r,b,...)$, where $b\leq {r\choose 2}+1$, by Macaulay's theorem we easily have that $$h_{i}\leq (...((({r\choose 2}+1)^{<2>})^{<3>})...)^{<i-1>}={r+i-2\choose i}+1.$$ \indent The desired upper-bound for $h_{e-i}$ immediately follows from the equality $h_{e-i}=h_i$.\\\indent 
In order to show that there exists a level $h$-vector of the form $(1,r,...,h_{e-i},...,r,2)$ for any integer $h_{e-i}\in [r,\min \lbrace {r-2+i\choose i}+1,{r-2+e-i\choose e-i}+1\rbrace ]$, we reason as in Theorem 3.5, i): consider the module $M^{'}=\langle F^{'},G^{'}\rangle $, where we choose $G^{'}=y_r^e$ and $F^{'}$ in variables $y_1,...,y_{r-1}$, the sum of powers of $h_{e-i}-1$ generic linear forms. Using Inverse Systems and Lemma 3.1, the result easily follows.\\\indent 
ii). Let $a>r$. Let $F\in M$ have $r_F$ first derivatives. The Gorenstein $h$-vector of $\langle F\rangle $, $(1,r_F,...,h_{e-i}^{F},...,r_F,1)$, of course has $h_{e-i}^{F}\leq {r_F-1+i\choose i}$. By Corollary 2.2, $(2,a,...,h_{e-i},...)-(1,r_F,...,h_{e-i}^{F},...)$ is an $O$-sequence. Therefore $h_{e-i}-h_{e-i}^{F}\leq {a-r_F-1+i\choose i}$, and the desired upper-bound for $h_{e-i}$ easily follows by Lemma 3.4.\\\indent 
In order to show that there exists a level $h$-vector of the form $(1,r,...,h_{e-i},...,a,2)$ for any choice of the integer $h_{e-i}$ in $[a, \min \lbrace {r-1+i\choose i}+{a-r-1+i\choose i},{r-1+e-i\choose e-i}\rbrace ]$, let us consider the Inverse System module $M^{'}=\langle F^{'},G^{'}\rangle $, where we choose $G^{'}$ only in variables $y_1,...,y_{a-r}$ and $F^{'}$ in all of the $r$ variables, both having degree $e$. The result follows by Lemmata 3.2 and 3.1 and the same kind of argument we used in proving Theorem 3.5, ii). This completes the proof of ii).\\\indent 
iii). If $a=r$ (now $r\leq 5$), choose $F^{'}$ only in variables $y_1,...,y_{r-1}$ as the sum of the $e$-th powers of $\alpha $ generic linear forms, where $\alpha \geq {r-2+\lfloor {e\over 2}\rfloor \choose \lfloor {e\over 2}\rfloor }$, and $G^{'}=y_r^e$. By Lemma 3.1, it is easy to check that the Inverse System module $M^{'}=\langle F^{'},G^{'}\rangle $ gives the $h$-vector $H$ where each entry is exactly the upper-bound of i).\\\indent 
If $a>r$ (for any $r\geq 2$), choose $G^{'}$ only in variables $y_1,...,y_{a-r}$ as the sum of powers of $\beta $ generic linear forms, where $\beta \geq {a-r-1+\lfloor {e\over 2}\rfloor \choose \lfloor {e\over 2}\rfloor }$, and $F^{'}$ in all of the $r$ variables as the sum of powers of $\gamma $ generic linear forms, where $\gamma \geq {r-1+\lfloor {e\over 2}\rfloor \choose \lfloor {e\over 2}\rfloor }$. By Lemmata 3.2 and 3.1, we can easily see that the Inverse System module $M^{'}=\langle F^{'},G^{'}\rangle $ gives the desired $h$-vector $H$.\\
\indent This concludes the proof of the theorem.{\ }{\ }\qed \\
\\\indent
{\bf Definition-Remark 3.9.} Consider the set of all the algebras with given codimension $r$ and socle-vector $s$ whose $h$-vector has an additional fixed entry, say $h_i$, with $1<i<e$. We define the algebras from this set having the (entry by entry) maximal $h$-vector, if this exists, as {\it quasi-compressed} (with respect to the data $(r,s,h_i)$).\\\indent
Therefore, Theorem 3.8 studies quasi-compressed algebras for level socle-vectors of type 2, codimension $r$, and for which the next to last entry $h_{e-1}=a$ of the $h$-vector is fixed.\\\indent
It would be interesting to continue and develop the study of quasi-compressed algebras that we have just begun, especially determining under which conditions they cannot exist.\\
\\\indent
{\bf Example 3.10.} In this example, we exhibit the maximum for the $h$-vectors $h=(1,5,h_2,h_3,h_4,h_5,h_6,5,2)$ of level algebras having codimension $r=5$, type 2, socle degree $e=8$, and the next to last entry $a$ of the $h$-vector equal to 5.\\\indent
By Theorem 3.8, the maximum for these $h$-vectors is $$H=(1,5,11,21,36,21,11,5,2).$$ \indent In order to construct a quasi-compressed level algebra with $h$-vector $H$, choose $F,G\in S=k[y_1,...,y_5]$ such that $G$ involves only $y_1,...,y_4$ and is the sum of the 8-th powers of (at least) 35 generic linear forms, and $F=y_5^8$. Then, by Lemma 3.1, the $R$-submodule $M=\langle F,G\rangle $ of $S$ gives the $h$-vector $H$ above.\\
\\\indent
{\bf Example 3.11.} We now exhibit the maximum for the $h$-vectors $h=(1,3,h_2,h_3,h_4,h_5,a,2)$ of level algebras having codimension $r=3$, type 2, socle degree $e=7$, and for which the next to last entry of the $h$-vector is $a$ (notice that, by Theorem 3.3, we have $3\leq a\leq 6$).\\\indent
Let $a=3$. By Theorem 3.8, the maximum for the level $h$-vectors above is $$H=(1,3,4,5,5,4,3,2).$$ \indent To construct a quasi-compressed level algebra with $h$-vector $H$, choose $F,G\in S=k[y_1,y_2,y_3]$ such that $G$ is only in variables $y_1,y_2$ and is the sum of the 7-th powers of (at least) 4 generic linear forms, and $F=y_3^7$. (See also Theorem 2.9, where the level $h$-vectors of the form $(1,3,...,3,2)$ are characterized.)\\\indent
Let $a=4$. By Theorem 3.8, the maximum for the level $h$-vectors above is $$H=(1,3,6,10,11,7,4,2).$$ \indent To construct a quasi-compressed level algebra with $h$-vector $H$, choose $F,G\in S$ such that $G=y_1^7$ and $F$ is in all of the three variables and is the sum of powers of (at least) 10 generic linear forms.\\\indent
Let $a=5$. By Theorem 3.8, the maximum for the level $h$-vectors above is $$H=(1,3,6,10,14,9,5,2).$$ \indent To construct a quasi-compressed level algebra with $h$-vector $H$, choose $F,G\in S$ such that $G$ is only in variables $y_1,y_2$ and is the sum of powers of (at least) 4 generic linear forms, and $F$ is in all of the three variables and is the sum of powers of (at least) 10 generic linear forms.\\\indent
Let $a=6$. By Theorem 3.8, the maximum for the level $h$-vectors above is $$H=(1,3,6,10,15,12,6,2).$$ \indent To construct a quasi-compressed level algebra with $h$-vector $H$, choose $F,G\in S$ such that both $F$ and $G$ are in all of the three variables and are the sum of powers of (at least) 8 generic linear forms each. In this case, since $a=6$ is the maximum possible value for the next to last entry $a$ of our level $h$-vectors, the level algebras having $h$-vector $H$ are actually {\it compressed} (see $[Ia2]$, $[FL]$, $[Za]$ and $[Za2]$).\\
\\
\section{The MFR's associated to a class of level $h$-vectors of codimension 3}
\indent 

In Theorem 2.9 we have characterized the level $h$-vectors of the form $h=(1,3,...,3,2)$: $h$ is level if and only if $h$ can be written as the sum of a Gorenstein $h$-vector of codimension 2 and $(0,1,1,...,1)$. In this section we characterize the possible minimal free resolutions associated to each of the $h$-vectors above.\\
\indent In order to do this, let us recall some basic facts. As above, let $R=k[x_1,...,x_r]$ and let $I$ be a homogeneous ideal of $R$ having no non-zero forms of degree 1. The {\it minimal free resolution} ({\it MFR}) of a standard graded algebra $A=R/I$ having {\it depth} 0 (i.e. all the homogeneous polynomials of $R$ are zero-divisors modulo $I$) is an exact sequence of $R$-modules of the form: $$0\longrightarrow F_r\longrightarrow F_{r-1}\longrightarrow ...\longrightarrow F_1\longrightarrow R\longrightarrow R/I \longrightarrow 0,$$
where, for $i=1,...,r$, $$F_i=\bigoplus_{j=1}^{n_i}R^{\beta_{i,j}}(-j),$$  and all the homomorphisms have degree 0.\\\indent The $\beta_{i,j}$'s are called the {\it graded Betti numbers} of $A$.\\\indent 
Then $\beta_{1,j}$ is the number of generators of $I$ in degree $j$. It is well-known that $F_r=\oplus_{j=1}^eR^{s_j}(-j-r)\neq 0,$ where $s(A)=s=(s_0,...,s_e)$ is, as usual, the socle-vector of $A$ (if $A$ is not artinian the socle is still defined as the annihilator of the maximal homogeneous ideal $\overline{m}=(\overline{x_1},...,\overline{x_r})$). Hence, the socle-vector may also be computed by considering the graded Betti numbers of the last module of the MFR. In particular, an artinian algebra $A$ is level of socle degree $e$ and type $s_e$ if and only if $F_r=R^{s_e}(-e-r)$.\\\indent 
For any algebra $A$, let $H(z)=\sum_{i=0}^{\infty }h_iz^i$ be the {\it Hilbert series} of $A$, where $h_i=\dim_k A_i$. ($A$ is artinian if and only if $h_i$ is eventually 0.) The growth condition of Macaulay's theorem, namely $h_{d+1}\leq h_d^{<d>}$ for $d\geq 1$, holds for the coefficients of the Hilbert series of any standard graded algebra (in Theorem 1.2 we have only stated that result for artinian algebras).\\\indent 
It can be shown that $H$ may be written as $H(z)={h(z)\over (1-z)^d}$, where $h$ is a polynomial having integral coefficients such that $h(1)\neq 0$, and $d$ is the {\it dimension} of $A$. In particular, $A$ is artinian if and only if $d=0$.\\\indent 
The MFR and the Hilbert series of $A$ are related by the following well-known formula (see, e.g., $[FL]$, p. 131, point (j) for a proof): \begin{equation}\label{sh}H(z)(1-z)^r=1+\sum_{i,j}(-1)^i\beta_{i,j}z^j.\end{equation}
\indent It follows from (\ref{sh}) and the observations above that, if $r=3$, determining the MFR of an artinian algebra $A=R/I$ with a given $h$-vector and socle-vector is equivalent to determining the degrees of the generators of $I$, i.e. only the first module of the MFR of $A$.\\
\\\indent 
Now fix the codimension $r=3$ and the socle-vector $(0,0,...,0,s_e=1)$, and let $T$ be a possible Gorenstein $h$-vector for this pair $(3,(0,0,...,0,1))$ (see Proposition 2.7 for Stanley's characterization). Let $${\bf F}: {\ }{\ }0\longrightarrow F_3\longrightarrow F_2\longrightarrow F_1\longrightarrow R\longrightarrow R/I \longrightarrow 0$$ be the MFR of a Gorenstein artinian algebra having data $(3,(0,0,...,0,1))$ as above and $h$-vector $T$. Then we can write $F_3=R(-e-3)$, $F_2=R(-p_1)\oplus ...\oplus R(-p_n)$ and $F_1=R(-q_1)\oplus ...\oplus R(-q_n)$, with $q_1\leq q_2\leq ...\leq q_n$ and $p_1\geq p_2\geq ...\geq p_n$.\\\indent 
For $i=1,...,n$, put $r_i=p_i-q_i$. Furthermore, let $k$ be the least degree in which $T$ is not {\it generic} (i.e. $h_k<{2+k\choose k}$), and let $\Delta^3T=(1,d_1,...,d_{e+3})$ be the third difference of $T$. Finally, put $\mu =2 \lceil {-\sum_{d_i<0}d_i\over 2}\rceil -1$. It is easy to see from (\ref{sh}) that, for each $i<e+3$ such that $d_i<0$, any Gorenstein algebra having $h$-vector $T$ must have at least $-d_i$ generators in degree $i$.\\\indent 
We are now ready to state the theorem of Diesel that characterizes the possible MFR's associated to the Gorenstein $h$-vectors of codimension 3.\\
\\\indent 
{\bf Theorem 4.1} ($[Di]$). {\it With the notation above, ${\bf F}$ is the minimal free resolution of some Gorenstein artinian algebra having data $(3,(0,0,...,0,1))$ and $h$-vector $T$ if and only if all of the following hold: the graded Betti numbers of ${\bf F}$ satisfy the functional equation $T(z)(1-z)^3=1+\sum_{i,j}(-1)^i\beta_{i,j}z^j$, $p_i+q_i=e+3$ for $i=1,...,n$, $n$ is odd, $\mu \leq n\leq 2k+1$, $r_1>0$, and $r_i+r_{n-i+2}>0$ for $i=2,...,{n+1\over 2}$.\\}
\\\indent 
{\bf Proof.} See $[Di]$. See also $[GPS]$.{\ }{\ }\qed \\
\\\indent 
We will also need Gotzmann's Persistence Theorem (Theorem 4.2 below). Recall that a vector space $V\subseteq R_d$ of forms of degree $d$ is called {\it Gotzmann} if $R_1V$ has the minimal possible dimension as a $k$-vector space given the dimension of $V$, i.e. if $$({r-1+d\choose d}-\dim_kV )^{<d>}={r-1+d+1\choose d+1}-\dim_kR_1V.$$
\\\indent 
{\bf Theorem 4.2} ($[Go]$). {\it Let $V$ be as above. If $V$ is Gotzmann, then $R_1V$ is also Gotzmann.\\}
\\
\indent {\bf Proof.} See $[Go]$, Satz 1 (where this result is shown in a more general context).{\ }{\ }\qed \\
\\\indent 
The following lemma (which is a result of M. Roth, personal communication) determines the level algebras of type 2, codimension 3 and socle degree $e$ having no Gorenstein quotient of codimension 3 and socle degree $e$. We omit its proof.\\
\\\indent 
{\bf Lemma 4.3} (Roth). {\it Let $M=\langle F,G\rangle \subseteq S=k[y_1,y_2,y_3]$ be the Inverse System module associated to a level algebra $A$ having $h$-vector $h=(1,3,...,3,2)$ and socle degree $e$. If no Gorenstein quotient of $A$ of socle degree $e$ has codimension 3, then, after a change of variables, we have $F=y_1y_3^{e-1}$ and $G=y_2y_3^{e-1}$. In particular, $h=(1,3,3,...,3,3,2)$.\\}
\\\indent 
Now we are ready for the main result of this section, that is the characterization of the MFR's associated to the level $h$-vectors of the form $h=(1,3,...,3,2)$ (we determined these $h$-vectors in Theorem 2.9). As we mentioned above, it is enough to determine the first module, $F_1$, of these MFR's ${\bf F}$.\\
\\\indent 
{\bf Theorem 4.4.} {\it Let $h=(1,3,...,3,2)$ be a level $h$-vector of socle-degree $e$, and let $j=h_{\lfloor e/2\rfloor }$. If $j>3$, then there is exactly one minimal free resolution ${\bf F}$ associated to $h$. It has $$F_1=R^2(-2)\bigoplus R(-(j-1))\bigoplus R(-(e-j+3))\bigoplus R(-(e+1)).$$\indent 
If $j=3$ (i.e. if $h=(1,3,3,...,3,3,2)$), then there are exactly two minimal free resolutions ${\bf F}$ associated to $h$. They have $$F_1=R^3(-2)\bigoplus R(-e)\bigoplus R^{\gamma }(-(e+1)),$$ with $\gamma =0$ or $\gamma =1$.\\}
\\\indent 
{\bf Proof.} Let $A=R/I$ be a level artinian algebra of socle degree $e$ having $h$-vector $h=(1,3,...,3,2)$, where $R=k[x_1,x_2,x_3]$, and let $M=\langle F,G\rangle \subseteq S=k[y_1,y_2,y_3]$ be the Inverse System module associated to $A$.\\
\indent Suppose first that no Gorenstein quotient of $A$ with socle degree $e$ has codimension 3. Then, by Lemma 4.3, after a change of variables, $M$ can be written as $M=\langle F,G\rangle $, where $F=y_1y_3^{e-1}$ and $G=y_2y_3^{e-1}$, and thus $h=(1,3,3,...,3,3,2)$.\\
\indent Notice that the annihilator of an Inverse System module generated by monomials is a monomial ideal, whose monomials in each degree are those of $R$ except for the monomials that (written in the $y_i$'s) belong to the module. From this observation, it is easy to see that the annihilator of $M$ above is $I=\langle x_1^2,x_1x_2,x_2^2,x_3^e\rangle $. This implies that the first module of the MFR of $R/I$ is that of the statement of the theorem for $j=3$ and $\gamma =0$.\\\indent 
Thus, from now on, we may suppose that $A$ has a Gorenstein quotient of socle degree $e$ and codimension 3; in other words, without loss of generality, we may suppose that the form $F$ has 3 linearly independent first derivatives. Hence, if $J\subseteq R$ is the Gorenstein ideal that annihilates $\langle F\rangle $, we have that $I_i=J_i$ for every $i\leq e-1$, and $J$ has exactly one generator more than $I$ in degree $e$. In particular, $J$ must have at least one generator in degree $e$.\\\indent 
By Theorem 2.9, the level $h$-vectors of the form $h=(1,3,...,3,2)$ are all and only those included in one of the following five cases:\\
\indent {\it Case 1)}. $h=(1,3,4,...,{e+4\over 2},...,4,3,2)$, where ${e+4\over 2}\geq 5$ is the maximal entry and occurs once (as $h_{e\over 2}$).\\
\indent {\it Case 2)}. $h=(1,3,4,...,{e+3\over 2},{e+3\over 2},...,4,3,2)$, where ${e+3\over 2}\geq 5$ is the maximal entry and occurs twice (the first time as $h_{e-1\over 2}$).\\
\indent {\it Case 3)}. $h=(1,3,4,...,t+2,t+2,...,t+2,...,4,3,2)$, where $t+2\geq 5$ is the maximal entry and occurs at least three times (the first as $h_t$).\\
\indent {\it Case 4)}. $h=(1,3,4,...,3,2)$, where the maximal entry is 4.\\
\indent {\it Case 5)}. $h=(1,3,...,3,2)$, where the maximal entry is 3.\\
\indent Let us consider {\it Case 1)}. The Gorenstein $h$-vector given by $\langle F\rangle $ is clearly $h^{'}=(1,3,4,...,{e+4\over 2},...,4,3,1)$. We have $$\Delta^3h^{'}=(1,0,-2,1,0,...,0,-2,2,0,...,0,-1,2,0,-1),$$ where the second \lq \lq -2" occurs in degree ${e+2\over 2}$.\\
\indent Using Diesel's Theorem 4.1 and $\Delta^3h^{'}$, we deduce that the number of generators of $J$ must be 5. By looking at the negative entries of $\Delta^3h^{'}$, we see that these generators occur in degrees $2,2,{e+2\over 2},{e+2\over 2},e$. Therefore the first module of the MFR of $R/J$ is $$R^2(-2)\bigoplus R^2(-{e+2\over 2})\bigoplus R(-e).$$ \indent Since $I$ and $J$ coincide in degrees less than $e$, we easily have that the first module of the MFR of $R/I$ must have the form $$F_1=R^2(-2)\bigoplus R^2(-{e+2\over 2})\bigoplus R^{\alpha }(-(e+1)),$$ for some non-negative integer $\alpha $.\\\indent 
{\it Claim.} $1\leq \alpha \leq 2$.\\\indent 
{\it Proof of claim.} We first show that $\alpha \geq 1$. From (\ref{sh}), it is easy to see that, in the MFR of $R/I$, $\beta_{2,e+1}-\beta_{1,e+1}=-3\cdot 2+3\cdot 3-1\cdot 4=-1$, whence $\beta_{1,e+1}=\alpha \geq 1.$\\\indent 
Let us now prove that $\alpha \leq 2$. Consider the Hilbert series $H^{''}$ of the algebra $A^{''}$ obtained from $A$ by getting rid of its $\alpha $ generators of degree $e+1$. We have $h^{''}_e=2$ and $h^{''}_{e+1}=\alpha $, whence $\alpha \leq 2$, in order to satisfy the growth condition of Macaulay's theorem. This proves the claim.\\\indent 
We will see later that $\alpha =2$ can never occur. This clearly completes the proof of this case, since here $j={e+4\over 2}$.\\\indent 
{\it Case 2.} The Gorenstein $h$-vector given by $\langle F\rangle $ is $h^{'}=(1,3,4,...,{e+3\over 2},{e+3\over 2},...,4,3,1)$. We have $$\Delta^3h^{'}=(1,0,-2,1,0,...,0,-1,0,1,0,...,0,-1,2,0,-1),$$ where the first \lq \lq -1" appears in degree ${e+1\over 2}$ and the second in degree $e$. By Diesel's theorem, $J$ has 5 generators, and only one of them is not indicated by the entries of $\Delta^3h^{'}$. But, by the formula $p_i+q_i=e+3$ of Diesel's theorem, we can easily see that the only degree possible for this last generator is ${e+3\over 2}$, since increasing by 1 any other graded Betti number $\beta_{1,j}$, $j\neq {e+3\over 2}$, would also require us to increase the Betti number $\beta_{1,e+3-j}$ by 1.\\\indent 
Hence the first module of the MFR of $R/J$ is $$R^2(-2)\bigoplus R(-{e+1\over 2})\bigoplus R(-{e+3\over 2})\bigoplus R(-e).$$ \indent Therefore the first module of the MFR of $R/I$ must have the form $$F_1=R^2(-2)\bigoplus R(-{e+1\over 2})\bigoplus R(-{e+3\over 2})\bigoplus R^{\alpha }(-(e+1)),$$ for some integer $\alpha $. The same argument as above shows that $1\leq \alpha \leq 2$, and we will see later that $\alpha \neq 2$. This completes the proof of this case, since $j={e+3\over 2}$.\\\indent 
{\it Case 3.} The Gorenstein $h$-vector given by $\langle F\rangle $ is $h^{'}=(1,3,4,...,t+2,t+2,...,t+2,...,4,3,1)$. We have $$\Delta^3h^{'}=(1,0,-2,0,1,0,...,0,-1,1,0,...,0,-1,1,0,...0,-1,2,0,-1),$$ where the first \lq \lq -1" appears in degree $t+1$ and the second in degree $e-t+1$. Thus, the number of generators of $J$ is again 5, and the only possible MFR for $R/J$ has the first module equal to $$R^2(-2)\bigoplus R(-(t+1))\bigoplus R(-(e-t+1))\bigoplus R(-e).$$ \indent Therefore the first module of the MFR of $R/I$ has the form $$F_1=R^2(-2)\bigoplus R(-(t+1))\bigoplus R(-(e-t+1))\bigoplus R^{\alpha }(-(e+1)),$$ for some integer $\alpha $. As above, $1\leq \alpha \leq 2$, and we will show later that $\alpha \neq 2$. This completes the proof of this case, since here $j=t+2$.\\\indent 
{\it Case 4.} Let us consider first the case $e\geq 6$. Then the Gorenstein $h$-vector given by $\langle F\rangle $ is $h^{'}=(1,3,4,4,...,4,3,1)$. We have $$\Delta^3h^{'}=(1,0,-2,0,1,0,...,0,-1,0,2,0,-1).$$ \indent From this formula for $\Delta^3h^{'}$, we deduce that $J$ might have either 3 or 5 generators (of which at least 2 in degree 2 and 1 in degree $e-1$), but, as we observed above, $J$ must also have a generator in degree $e$. Hence the number of generators of $J$ can only be 5. We need to find the last generator. Suppose it is in degree $e+3-i$. Then the Gorenstein MFR has its first module equal to $$R^2(-2)\bigoplus R(-(e+3-i))\bigoplus R(-(e-1))\bigoplus R(-e)$$ and the second module equal to $$R(-3)\bigoplus R(-4)\bigoplus R(-i)\bigoplus R^2(-(e+1)),$$ where $i\geq 4$. In fact, if $i=3$, the inequality $r_3+r_4>0$ of Diesel's theorem is: $4-(e-1)+3-e>0$, i.e. $e<4$, which is impossible. Hence $i\geq 4$.\\\indent 
Therefore $r_3+r_4>0$ is: $i-(e+3-i)+4-(e-1)>0$, i.e. $i>e-1$. Since $i<e+1$ ($J$ has exactly 2 generators in degree 2), we have $i=e$. Thus the last generator of $J$ is in degree 3.\\\indent 
Hence the MFR of $R/J$ has the first module equal to $$R^2(-2)\bigoplus R(-3)\bigoplus R(-(e-1))\bigoplus R(-e).$$\indent 
Therefore the first module of the MFR of $R/I$ has the form $$F_1=R^2(-2)\bigoplus R(-3)\bigoplus R(-(e-1))\bigoplus R^{\alpha }(-(e+1)),$$ for some integer $\alpha $. As above, $1\leq \alpha \leq 2$, and we will show later that $\alpha \neq 2$.\\\indent
Now let $e=4$. Thus $h^{'}=(1,3,4,3,1)$, and we have $$\Delta^3h^{'}=(1,0,-2,-1,1,2,0,-1).$$\indent The argument we gave above for $e\geq 6$ also holds for $e=4$, and therefore we are done.\\\indent
Finally, let $e=5$. Hence $h^{'}=(1,3,4,4,3,1)$, and we have $$\Delta^3h^{'}=(1,0,-2,0,0,0,2,0,-1).$$\indent The same reasoning as above, since $J$ must have at least one generator in degree 5, shows that $J$ has exactly 5 generators. The formula for $\Delta^3h^{'}$ only indicates that $J$ has 2 generators in degree 2. Hence we have to find the last two generators. Suppose they are in degrees $l$ and $m$, with $3\leq l\leq m\leq 5$. Writing down the MFR for $R/J$, we immediately see that the inequality $r_3+r_4>0$ of Diesel's theorem is: $(8-l-l)+(8-m-m)>0$. Hence, we can only have $l=3$ and $m=4$, or $l=m=3$.\\\indent
However, if $l=m=3$, it is easy to see that equation (\ref{sh}) is not satisfied, and therefore the only possible choice for the degrees of the last two generators of $J$ is $l=3$ and $m=4$. This leads us to the same formulas for the MFR's of $R/J$ and $R/I$ that we gave above for $e\geq 6$. Thus (up to showing, as usual, that $\alpha \neq 2$), we have concluded the argument for $e=5$.\\\indent
This clearly completes the proof of Case 4, since $j=4$.\\\indent 
{\it Case 5.} Consider first the case $e\geq 4$. The Gorenstein $h$-vector given by $\langle F \rangle $ is $h^{'}=(1,3,3,...,3,1)$. Therefore we have $$\Delta^3h^{'}=(1,0,-3,2,0,...,0,-2,3,0,-1).$$ \indent Hence, by Diesel's theorem, $J$ has 5 generators and the first module of the MFR of $R/J$ is $$R^3(-2)\bigoplus R^2(-e).$$ \indent Therefore the MFR of $R/I$ has the first module equal to $$F_1=R^3(-2)\bigoplus R(-e)\bigoplus R^{\gamma }(-(e+1)),$$ for some non-negative integer $\gamma $. The same argument we used above for $\alpha $ also shows that $\gamma \leq 2$ (but cannot show that $\gamma \neq 0$). We will see next that $\gamma =0$ and $\gamma =1$ may both occur in the first module of the MFR of a level algebra having $h$-vector $h=(1,3,3,...,3,3,2)$. Instead, we will show that $\gamma \neq 2$.\\\indent
Now let $e=3$, i.e. $h^{'}=(1,3,3,1)$. Thus $$\Delta^3h^{'}=(1,0,-3,0,3,0,-1).$$ \indent
Reasoning as in the proof of Case 4, since $J$ must have at least one generator in degree $e=3$, we have that $J$ has exactly 5 generators. Of course three of these generators must be in degree 2 and two in degree 3. The rest of the argument does not differ from that for $e\geq 4$ and will be omitted.\\\indent
This clearly concludes the proof of Case 5 and that of the theorem, up to showing the next two lemmata:\\
\\\indent 
{\bf Lemma 4.5.} {\it $\gamma =0$ and $\gamma =1$ may both occur.\\}
\\\indent 
{\bf Lemma 4.6.} {\it $\gamma \neq 2$ and $\alpha \neq 2$.\\}
\\\indent 
{\bf Proof of Lemma 4.5.} The fact that $\gamma =0$ may occur has been shown at the beginning of the proof of Theorem 4.4, when we considered the case of a level algebra $A$ having no Gorenstein quotients of socle degree $e$ and codimension 3.\\
\indent Let us now show that $\gamma =1$ may also occur. Consider the Inverse System module $M=\langle y_1^{e-1}y_2,y_3^e\rangle \subseteq S$. Then the algebra $R/I$, where $I$ is the annihilator of $M$, has $h$-vector (1,3,3,...,3,3,2). Moreover, by the considerations we made earlier about the annihilator of a module generated by monomials, we can easily deduce that $I=\langle x_1x_3,x_2x_3,x_2^2,x_1^e,x_3^{e+1}\rangle $. Therefore the first module of the MFR of $R/I$ has the Betti number $\gamma =1$. This completes the proof of the lemma.{\ }{\ }\qed \\
\\\indent 
{\bf Proof of Lemma 4.6.} We want to show, with the notation above, that $\alpha =2$ may never occur (the proof of $\gamma \neq 2$ requires no modifications). Suppose that $\alpha =2$ occurs for some level algebra $A$. Consider the algebra $A^{''}$ obtained from $A$ by getting rid of the 2 generators of $A$ of degree $e+1$. Then the coefficients of the Hilbert series of $A^{''}$ are $1,3,...,3,2,2,...$, where \lq \lq 2" is repeated from degree $e$ on. In fact, $A$ and $A^{''}$ coincide up to degree $e$, and afterwards $A^{''}$ has no more generators. Thus, since $2^{<i>}=2$ for every $i\geq e$, by Gotzmann's Persistence Theorem 4.2 all the entries from degree $e$ on must be equal to 2.\\\indent 
Standard considerations show that the algebra $A^{''}$ has dimension 1 and depth 0. Hence, as we observed at the beginning of this section, the last module of the MFR of $A^{''}$ represents its socle.\\\indent 
{\it Claim.} In the MFR of $A^{''}$, $\beta_{3,e+2}>0$.\\\indent 
{\it Proof of claim.} An easy computation shows that, in the l.h.s. of formula (\ref{sh}) for $A^{''}$, the coefficient of $z^{e+2}$ is -1. Hence, since $A^{''}$ has no generators of degree $e+2$, we have $\beta_{2,e+2}-\beta_{3,e+2}=-1$, and therefore $\beta_{3,e+2}>0$. This proves the claim.\\\indent 
Hence $A^{''}$ has a non-zero socle in degree $e+2-3=e-1$. But $A$ and $A^{''}$ coincide in degrees less than or equal to $e$, and therefore $A$ has also a non-zero socle in degree $e-1$, a contradiction since $A$ is level. This concludes the proof of the lemma and that of the theorem.{\ }{\ }\qed \\
\\\indent 
{\bf Acknowledgements.} We wish to express our warm gratitude to Professor Anthony Iarrobino for sending us the first draft of $[Ia]$ and for his interesting comments on a previous version of this work, and to Professor Mike Roth for giving us a copy of his result (Lemma 4.3).\\
\\
\\
\\
{\bf \huge References}\\
\\
$[BI]$ {\ } D. Bernstein and A. Iarrobino: {\it A non-unimodal graded Gorenstein Artin algebra in codimension five}, Comm. in Algebra 20 (1992), No. 8, 2323-2336.\\
$[BG]$ {\ } A.M. Bigatti and A.V. Geramita: {\it Level Algebras, Lex Segments and Minimal Hilbert Functions}, Comm. in Algebra 31 (2003), 1427-1451.\\
$[Bo]$ {\ } M. Boij: {\it Graded Gorenstein Artin algebras whose Hilbert functions have a large number of valleys}, Comm. in Algebra 23 (1995), No. 1, 97-103.\\
$[BL]$ {\ } M. Boij and D. Laksov: {\it Nonunimodality of graded Gorenstein Artin algebras}, Proc. Amer. Math. Soc. 120 (1994), 1083-1092.\\
$[BH]$ {\ } W. Bruns and J. Herzog: {\it Cohen-Macaulay rings}, Cambridge studies in advanced mathematics, No. 39, Revised edition (1998), Cambridge, U.K..\\
$[CI]$ {\ } Y.H. Cho and A. Iarrobino: {\it Hilbert Functions and Level Algebras}, J. of Algebra 241 (2001), 745-758.\\
$[CI2]$ {\ } Y.H. Cho and A. Iarrobino: {\it Inverse Systems of Zero-dimensional Schemes in ${\bf P}^n$}, J. of Algebra, to appear.\\
$[Di]$ {\ } S.J. Diesel: {\it Irreducibility and dimension theorems for families of height 3 Gorenstein algebras}, Pac. J. Math. 172 (1996), No. 2, 365-397.\\
$[FL]$ {\ } R. Fr\"oberg and D. Laksov: {\it Compressed Algebras}, Conference on Complete Intersections in Acireale, Lecture Notes in Mathematics, No. 1092 (1984), 121-151, Springer-Verlag.\\
$[Ge]$ {\ } A.V. Geramita: {\it Inverse Systems of Fat Points: Waring's Problem, Secant Varieties and Veronese Varieties and Parametric Spaces of Gorenstein Ideals}, Queen's Papers in Pure and Applied Mathematics, No. 102, The Curves Seminar at Queen's (1996), Vol. X, 3-114.\\
$[GHMS]$ {\ } A.V. Geramita, T. Harima, J. Migliore and Y.S. Shin: {\it The Hilbert Function of a Level Algebra}, Memoir of the Amer. Math. Soc., to appear.\\
$[GHS]$ {\ } A.V. Geramita, T. Harima and Y.S. Shin: {\it Some special configurations of points in ${\bf P}^n$}, J. of Algebra 268 (2003), No. 2, 484-518.\\
$[GPS]$ {\ } A.V. Geramita, M. Pucci and Y.S. Shin: {\it Smooth Points of $GOR(T)$}, J. of Pure and Applied Algebra 122 (1997), 209-241.\\
$[Go]$ {\ } G. Gotzmann: {\it Eine Bedingung f\"ur die Flachheit und das Hilbertpolynom cines graduierten Ringes}, Math. Z. 158 (1978), No. 1, 61-70.\\
$[Ha]$ {\ } T. Harima: {\it Some examples of unimodal Gorenstein sequences}, J. of Pure and Applied Algebra 103 (1995), No. 3, 313-324.\\
$[Ia]$ {\ } A. Iarrobino: {\it Hilbert functions of Gorenstein algebras associated to a pencil of forms}, preprint (math.AC/0412361); Proc. of the Conference on Projective Varieties with Unexpected Properties (Siena 2004), C. Ciliberto {\it et al.} eds. (de Gruyter), to appear.\\
$[Ia2]$ {\ } A. Iarrobino: {\it Compressed Algebras: Artin algebras having given socle degrees and maximal length}, Trans. Amer. Math. Soc. 285 (1984), 337-378.\\
$[Ia3]$ {\ } A. Iarrobino: {\it Ancestor ideals of vector spaces of forms, and level algebras}, J. of Algebra 272 (2004), 530-580.\\
$[IK]$ {\ } A. Iarrobino and V. Kanev: {\it Power Sums, Gorenstein Algebras, and Determinantal Loci}, Springer Lecture Notes in Mathematics (1999), No. 1721, Springer, Heidelberg.\\
$[IS]$ {\ } A. Iarrobino and H. Srinivasan: {\it Some Gorenstein Artin algebras of embedding dimension four, I: components of $PGOR(H)$ for $H=(1,4,7,...,1)$}, J. of Pure and Applied Algebra, to appear.\\
$[Ma]$ {\ } F.H.S. Macaulay: {\it The Algebraic Theory of Modular Systems}, Cambridge Univ. Press, Cambridge, U.K. (1916).\\
$[MN]$ {\ } J. Migliore and U. Nagel: {\it Reduced arithmetically Gorenstein schemes and simplicial polytopes with maximal Betti numbers}, Adv. Math. 180 (2003), 1-63.\\
$[So]$ {\ } J. S\"oderberg: {\it Artinian Level Modules and Cancellable Sequences}, J. of Algebra 280 (2004), No. 2, 610-623.\\
$[St]$ {\ } R. Stanley: {\it Hilbert functions of graded algebras}, Adv. Math. 28 (1978), 57-83.\\
$[St2]$ {\ } R. Stanley: {\it Cohen-Macaulay Complexes, Higher Combinatorics}, M. Aigner Ed., Reidel, Dordrecht and Boston (1977), 51-62.\\
$[Za]$ {\ } F. Zanello: {\it Extending the idea of compressed algebra to arbitrary socle-vectors}, J. of Algebra 270 (2003), No. 1, 181-198.\\
$[Za2]$ {\ } F. Zanello: {\it Extending the idea of compressed algebra to arbitrary socle-vectors, II: cases of non-existence}, J. of Algebra 275 (2004), No. 2, 730-748.

}

\end{document}